# Combining domain knowledge and statistical models in time series analysis

Tze Leung Lai[1] and Samuel Po-Shing Wong[2]

*Stanford University and The Chinese Universty of Hong Kong*

**Abstract:** This paper describes a new approach to time series modeling that combines subject-matter knowledge of the system dynamics with statistical techniques in time series analysis and regression. Applications to American option pricing and the Canadian lynx data are given to illustrate this approach.

## 1. Introduction

In their Fisher Lectures at the Joint Statistical Meetings, Cox [11] and Lehmann [31] mentioned two major types of stochastic models in statistical analysis, namely, *empirical* and *substantive* (or mechanistic). Whereas substantive models are explanatory and related to subject-matter theory on the mechanisms generating the observed data, empirical models are interpolatory and aim to represent the observed data as a realization of a statistical model chosen largely for its flexibility, tractability and interpretability but not on the basis of subject-matter knowledge. Cox [11] also mentioned a third type of stochastic models, called *indirect* models, that are used to evaluate statistical procedures or to suggest methods for analyzing complex data (such as hidden Markov models in image analysis). He noted, however, that the distinctions between the different types of models are important mostly when formulating and checking them but that these types are not rigidly defined, since "quite often parts of the model, e.g., those representing systematic variation, are based on substantive considerations with other parts more empirical." In this paper, we elaborate further the complementary roles of empirical and substantive models in time series analysis and describe a basis function approach to combining subject-matter (domain) knowledge with statistical modeling techniques.

This basis function approach was first developed in [29] for the valuation of American options. In Sections 2 and 3 we review the statistical and subject-matter models for option pricing in the literature as examples of empirical and substantive models in time series analysis. Section 4 describes a combined substantive-empirical approach via basis functions, in which the substantive component is associated with basis functions of a certain form, and the empirical component uses flexible and computationally convenient basis functions such as regression splines. The work of Lai and Wong [29] on option pricing and recent related work in financial time series are reviewed to illustrate this approach. Section 5 applies this approach to a widely studied data set in the nonlinear time series literature, namely, the Canadian lynx data set that records the annual numbers of Canadian lynx trapped in the Mackenzie River district from 1821 to 1934. We use substantive models from the ecology literature together with multivariate adaptive regression splines to come up with a new time series model for these data. Some concluding remarks are given in Section 6.

[1]Department of Statistics, Stanford Univeristy, Stanford, CA 94305, U.S.A., e-mail: lait@stat.stanford.edu
[2]Department of Statistics, The Chinese University of Hong Kong, Shatin, N.T., Hong Kong, e-mail: samwong@sta.cuhk.edu.hk

## 2. Statistical (empirical) time series models

The development of statistical time series models in the past fifty years has witnessed a remarkable confluence of basic ideas from various areas in statistics and probability, coupled with the powerful influence from diverse fields of applications ranging from economics and finance to signal processing and control systems. The first phase of this development was concerned with stationary time series, leading to MA (moving average), AR (autoregressive) and ARMA representations in the time domain and transfer function representations in the frequency domain. This was followed by extensions to nonstationary time series, either by fitting (not necessarily stationary) ARMA models or by the Box-Jenkins approach involving the ARIMA (autoregressive integrated moving average) models and their seasonal SARIMA counterparts. More general fractional differencing then led to the ARFIMA models. The next phase of the development was concerned with nonlinear time series models, beginning with bilinear models that add cross-product terms $y_{t-i}\epsilon_{t-j}$ to the usual ARMA model $y_t = \beta_1 y_{t-1} + \cdots + \beta_p y_{t-p} + \epsilon_t + c_1 \epsilon_{t-1} + \cdots + c_q \epsilon_{t-q}$, and threshold autoregressive and regime switching models that introduce nonlinearities into the usual autoregressive models via state-dependent changes or Markov jumps in the autoregressive parameters. The monograph by Tong [44] summarized these and other nonlinear time series models in the previous literature. The appropriateness of the parametric forms assumed in these nonlinear time series models, however, may be difficult to justify in real applications, as pointed out by Chen and Tsay [9].

Whereas the AR model $y_t = \beta_1 y_{t-1} + \cdots + \beta_p y_{t-p} + \epsilon_t$ is related to linear regression since $\boldsymbol{\beta}^T \mathbf{x}_t$ is the regression function $E(y_t|\mathbf{x}_t)$ of $y_t$ given $\mathbf{x}_t := (y_{t-1}, \ldots, y_{t-p})^T$, and likewise its nonlinear parametric extensions $y_t = f(\mathbf{x}_t, \boldsymbol{\beta}) + \epsilon_t$ are related to nonlinear regression, Chen and Tsay [9, 10] proposed to use nonparametric regression for $E(y_t|\mathbf{x}_t)$ instead. They started with functional-coefficient autoregressive (FAR) models of the form $y_t = f_1(\mathbf{x}_t^*)y_{t-1} + \cdots + f_p(\mathbf{x}_t^*)y_{t-p} + \epsilon_t$, where $f_1, \ldots, f_p$ are unspecified functions to be estimated by local linear regression and $\mathbf{x}_t^* = (y_{t-i_1}, \ldots, y_{t-i_d})^T$ with $i_1 < \cdots < i_d$ chosen from $\{1, \ldots, p\}$. Because of sparse data in high dimensions, local linear regression typically require $d$ to be 1 or 2. To deal with nonparametric regression in higher dimensions, they considered additive autoregressive models of the form $y_t = f_1(y_{t-i_1}) + \cdots + f_d(y_{t-i_d}) + \epsilon_t$, in which the $f_i$ can be estimated nonparametrically via the generalized additive model (GAM) of Hastie and Tibshirani [19] . Making use of Friedman's [15] multivariate adaptive splines (MARS), Lewis and Stevens [34] and Lewis and Ray [32, 33] developed spline models for empirical modeling of time series data. Weigend, Rummelhart and Huberman [48] and Weigend and Gershenfeld [47] proposed to use neural networks (NN) to model $E(y_t|\mathbf{x}_t)$, while Lai and Wong [28] considered a variant called stochastic neural networks, for which they could use the EM algorithm to develop efficient estimation procedures that have much lower computational complexity than those for conventional neural networks.

The preceding time series models are autonomous, relating the dynamics of $y_t$ to the past states. In econometrics and engineering, the outputs $y_t$ are related not only to the past outputs but also to the past inputs $u_{t-d}, \ldots, u_{t-k}$. Therefore the AR model has been extended to the ARX model (where X stands for exogenous inputs) $y_t = \boldsymbol{\beta}^T \mathbf{x}_t + \epsilon_t$ with $\mathbf{x}_t = (y_{t-1}, \ldots, y_{t-p}, u_{t-d}, \ldots, u_{t-k})^T$. Instead of assuming a linear or nonlinear parametric regression model, one can use nonparametric regression to estimate $E(y_t|\mathbf{x}_t)$, as in the following financial application.

**Example 1.** As noted by Ross [40], option pricing theory is "the most successful theory not only in finance, but in all of economics." A call (put) option gives the holder the right to buy (sell) the underlying asset (e.g. stock) by a certain date $T$ (known as the "expiration date" or "maturity") at a certain price (known as the "strike price" and denoted by $K$). European options can be exercised only on the expiration date, whereas American options can be exercised at any time up to the expiration date. The celebrated Black-Scholes theory, which will be reviewed in Section 3, yields the following pricing formulas for the prices $c_t$ and $p_t$ of European call and put options at time $t \in [0, T)$:

(2.1) $c_t = S_t e^{-d(T-t)} \Phi(d_1(S_t, K, T-t)) - K e^{-r(T-t)} \Phi(d_2(S_t, K, T-t))$,
(2.2) $p_t = K e^{-r(T-t)} \Phi(-d_2(S_t, K, T-t)) - S_t e^{-d(T-t)} \Phi(-d_1(S_t, K, T-t))$,

where $\Phi$ is the cumulative distribution function of the standard normal random variable, $S_t$ is the price of the underlying asset at time $t$, $d$ is the dividend rate of the underlying asset, $d_1(x, y, v) = \{\log(x/y) + (r - d + \sigma^2/2)v\}/\sigma\sqrt{v}$ and $d_2(x, y, v) = d_1(x, y, v) - \sigma\sqrt{v}$. Hutchinson, Lo and Poggio [22] pointed out that the success of the formulas (2.1) and (2.2) depends heavily on the specification of the dynamics of $S_t$. Instead of using any particular model of $S_t$, they proposed a data-driven way for pricing and hedging with a minimal assumption: independent increments of the underlying asset price. Noting that $y_t$ ($= c_t$ or $p_t$) is function of $S_t/K$ and $T - t$ with $r$ and $\sigma$ being constant, they assume $y_t = Kf(S_t/K, T - t)$ and approximate $f$ by taking $\mathbf{x}_t = (S_t/K, T - t)^T$ in the following models:

  (i) radial basis function (RBF) networks $f(\mathbf{x}) = \beta_0 + \boldsymbol{\alpha}^T \mathbf{x} + \sum_{i=1}^{I} \beta_i h_i(\|A(\mathbf{x} - \boldsymbol{\gamma}_i)\|)$, where $A$ is a positive definite matrix and $h_i$ is of the RBF type $e^{-u^2/\sigma_i^2}$ or $(u^2 + \sigma_i^2)^{1/2}$;
  (ii) neural networks $f(\mathbf{x}) = \psi(\beta_0 + \sum_{i=1}^{I} \beta_i h(\boldsymbol{\gamma}_i + \boldsymbol{\alpha}_i^T \mathbf{x}))$, where $h(u) = 1/(1 + e^{-u})$ is the logistic function and $\psi$ is either the identity function or the logistic function;
  (iii) projection pursuit regression (PPR) networks $f(\mathbf{x}) = \beta_0 + \sum_{i=1}^{I} \beta_i h_i(\boldsymbol{\alpha}_i^T \mathbf{x})$, where $h_i$ is an unspecified function that is estimated from the data by PPR.

The $\boldsymbol{\alpha}_i$, $\beta_i$ and $\boldsymbol{\gamma}_i$ above are unknown parameters of the network that are to be estimated from the data. As pointed out in [22], all three classes of networks have some form of "universal approximation property" which means their approximation bounds do not depend on the dimensionality of the predictor variable $x$; see [2]. It should be noted that the above transformation of $S_t$ to $S_t/K$ can be motivated not only from the assumption on $S_t$ but also from the special feature of options data. Although the strike price $K$ could be any positive number theoretically, the options exchange only sets strike prices at a multiple of a fundamental unit. For example, Chicago Board Options Exchange (CBOE) sets strike prices at multiples of $5 for stock prices in the $25 to $200 range. Also, only those options with strike prices closet to the current stock price are traded and thus their prices are observed. Since $S_t$ is non-stationary in general, the observed $K$ is also non-stationary. Such features

create sparsity of data in the space of $(S_t, K, T-t)$. Training the options pricing formula in the form of $f(S_t, K, T-t)$ can only interpolate the data and can hardly produce any good prediction because $(S_t, K)$ in the future can be very different from the data used in estimating $f$. The proposed transformation makes use of the fact that all observed and future $S_t/K$ are close to 1. Therefore, the proposed transformation captures the stationary structure of the data and enable the non-parametric models to predict well. Another point that Hutchinson, Lo and Poggio [22] highlighted is the measure of performance of the estimated pricing formula. According to their simulation study, even a linear $f(S_t/K, T-t)$ can give $R^2 \approx 90\%$ (Table I of Hutchinson, Lo and Poggio [22]). However, such a linear $f$ implies a constant delta hedging scheme which would provide poor hedging results. Since the primary function of options is hedging the risk created by changes in the price of the underlying asset, Hutchinson, Lo and Poggio [22] suggested using, instead of $R^2$, the hedging error measures $\xi = e^{-rT}E[|V(T)|]$ and $\eta = e^{-rT}[EV^2(T)]^{1/2}$, where $V(T)$ is the value of the hedged portfolio at time $T$. In a perfect Black-Scholes world, $V(T)$ should be 0 if Black-Scholes formula is used. However, from the simulation study, the Black-Scholes formulas still give $\xi > 0$ and $\eta > 0$ because time is discrete. Hutchinson, Lo and Poggio [22] reported that RBF, NN and PPR all give hedging measures comparable to those of the Black-Scholes in the simulation study. For real data analysis of futures options, RBF, NN and PPR performed better than the Black-Scholes formula in terms of hedging.

For American options, instead of using these learning networks to approximate the option price, Broadie et al. [5] used kernel smoothers to estimate the option pricing formula of an American option. Using a training sample of daily closing prices of American calls on the S&P100 Index that were traded on the Chicago Board Options Exchange from 3 January 1984 to 30 March 1990, they compared the nonparametric estimates of American call option prices at a set of $(S/K, t^*)$ values with corresponding parametric estimates obtained by using the approximations to American option prices due to Broadie and Detemple [4], and found significant differences between the parametric and nonparametric estimates.

## 3. Substantive (mechanistic) models

In control engineering, the dynamics of linear input-output systems are often given by ordinary differential equations, whose discrete-time approximations in the presence of noise have led to the ARX models (for white noise), and ARMAX models (for colored noise) in the preceding section. The problem of choosing the inputs sequentially so that the outputs are as close as possible to some target values when the model parameters are unknown and have to be estimated on-line has a large literature under the rubric of *stochastic adaptive control*; see Goodwin, Ramadge and Caines [16], Lai and Wei [27], Lai and Ying [30] and Guo and Chen [17]. More general dynamics in the presence of additive noise have led to stochastic differential equations (SDEs), whose discrete-time approximations are related to nonlinear time series models described in the preceding section. One such SDE is geometric Brownian motion (GBM) for the asset price process in the Black-Scholes option pricing theory. In view of Ito's formula, the GBM dynamics for the asset price $S_t$ translate into SDE dynamics for the option price $f(t, S_t)$. Such implied dynamics from the mechanistic model can be combined with subject-matter theory to derive the functional form or differential equation for $f$ and other important corollaries of the theory, as illustrated in the following.

**Example 2.** In the Black-Scholes model, the asset price $S_t$ is assumed to be GBM defined by the SDE

(3.1) $$dS_t/S_t = \mu dt + \sigma dw_t,$$

where $w_t, t \geq 0$, is Brownian motion. Letting $f(t, S)$ be the price of the option at time $t$ when $S_t = S$, it follows from (3.1) and Ito's formula that

$$\begin{aligned} df(t, S_t) &= \frac{\partial f}{\partial t}dt + \frac{\partial f}{\partial S}dS_t + \frac{1}{2}\frac{\partial^2 f}{\partial S^2}\sigma^2 S_t^2 dt \\ &= \left(\frac{\partial f}{\partial t} + \mu S_t \frac{\partial f}{\partial S} + \frac{1}{2}\sigma^2 S_t^2 \frac{\partial^2 f}{\partial S^2}\right)dt + \sigma S_t \frac{\partial f}{\partial S}dw_t. \end{aligned}$$

For simplicity assume that the asset does not pay dividends, i.e., $d = 0$. Consider an option writer's portfolio at time $t$, consisting of $-1$ option and $y_t$ units of the asset. The value of the portfolio $\pi_t$ is $-f(t, S_t) + y_t S_t$ and therefore

$$d\pi_t = -\left(\frac{\partial f}{\partial t} + \mu S_t \frac{\partial f}{\partial S} + \frac{1}{2}\sigma^2 S_t^2 \frac{\partial^2 f}{\partial S} - \mu y_t S_t\right)dt + \sigma S_t \left(y_t - \frac{\partial f}{\partial S}\right)dw_t.$$

Hence setting $y_t = \partial f/\partial S$ yields a risk-free portfolio. This is the basis of *delta hedging* in the options theory of Black and Scholes [3], who denote $\partial f/\partial S$ by $\Delta$. Besides GBM dynamics for the asset price, the Black-Scholes theory also assumes that there are no transaction costs and no limits on short selling and that trading can take place continuously so that delta hedging is feasible. Since economic theory prescribes absence of arbitrage opportunities in equilibrium, $\pi_t$ that consists of $-1$ option and $\Delta$ units of the asset should have the same return as $r\pi_t dt = r(-f + S_t \Delta)dt$, yielding the Black-Scholes PDE for $f$:

(3.2) $$\frac{\partial f}{\partial t} + rS\frac{\partial f}{\partial S} + \frac{1}{2}\sigma^2 S^2 \frac{\partial^2 f}{\partial S^2} = rf, \quad 0 \leq t < T,$$

with the boundary condition $f(T, S) = g(S)$, where $g(S) = (K - S)_+$ for a put option, and $g(S) = (S - K)_+$ for a call option, where $x_+ = \max(x, 0)$. This PDE has the explicit solution (2.1) or (2.2) with $d = 0$. If the asset pays dividend at rate $d$, then a modification of the preceding argument yields (3.2) in which $rS(\partial f/\partial S)$ is replaced by $(r - d)S(\partial f/\partial S)$.

Merton [37] extended the Black-Scholes theory for pricing European options to American options that can be exercised at any time prior to the expiration date. Optimal exercise of the option is shown to occur when the asset price exceeds (or falls below) an exercise boundary $\partial \mathcal{C}$ for a call (or put) option. The Black-Scholes PDE still holds in the continuation region $\mathcal{C}$ of $(t, S_t)$ before exercise, and $\partial \mathcal{C}$ is determined by the free boundary condition $\partial f/\partial S = 1$ (or $-1$) for a call (or put) option. Unlike the explicit formula (2.1) or (2.2) for European options, there is no closed-form solution of the free-boundary PDE and numerical methods such as finite differences are needed to compute American option prices under this theory.

By the Feynman-Kac formula, the PDE (3.2) has a probabilistic representation $f(t, S) = E[e^{-r(T-t)}g(S_T)|S_t = S]$, and the expectation $E$ is with respect to the "equivalent martingale measure" under which $dS_t/S_t = (r - d)dt + \sigma dw_t$. This representation generalizes to American options as the value function of the optimal stopping problem

(3.3) $$f(t, S) = \sup_{\tau \in \mathcal{T}_{t,T}} E[e^{-r(\tau-t)}g(S_\tau)|S_t = S]$$

where $\mathcal{T}_{t,T}$ denotes the set of stopping times $\tau$ taking values between $t$ and $T$. Cox, Ross and Rubinstein [12] proposed to approximate GBM by a binomial tree, with root node $S_0$ at time 0, so that (3.3) can be approximated by a discrete-time and discrete-state optimal stopping problem that can be solved by backward induction. Denote $f(t,S)$ by $C(t,S)$ for an American call option, and by $P(t,S)$ for an American put option. Jacka [23] and Carr, Jarrow and Myneni [7] derived the decomposition formula

$$
\begin{aligned}
P(t,S) = p(t,S) + K\rho e^{\rho u} \int_u^0 \Big\{ & e^{-\rho s}\Phi\Big(\frac{\bar{z}(s)-z}{\sqrt{s-u}}\Big) \\
- \theta e^{-(\theta\rho s+u/2)+z}\Phi\Big(& \frac{\bar{z}(s)-z}{\sqrt{s-u}} - \sqrt{s-u}\Big) \Big\}ds,
\end{aligned}
\tag{3.4}
$$

and a similar formula relating $C(t,S)$ to $c(t,S)$, where $\bar{z}(u)$ is the early exercise boundary $\partial\mathcal{C}$ under the transformation

$$
\rho = r/\sigma^2,\ \theta = d/r;\ u = \sigma^2(t-T),\ z = \log(S/K) - (\rho - \theta\rho - 1/2)u.
\tag{3.5}
$$

Ju [24] found that the early exercise premium can be computed in closed form if $\partial\mathcal{C}$ is a piecewise exponential function which corresponds to a piecewise linear $\bar{z}(u)$. By using such assumption, Ju [24] reported numerical studies showing his method with 3 equally spaced pieces substantially improves previous approximations to option prices in both accuracy and speed. AitSahlia and Lai [1] introduced the transformation (3.5) to reduce GBM to Brownian motion and showed that $\bar{z}(u)$ is indeed well approximated by a piecewise linear function with a few pieces. The integral obtained by differentiating that in (3.4) with respect to $S$ also has a closed-form expression when $\bar{z}(\cdot)$ is piecewise linear, and approximating $\bar{z}(\cdot)$ by a linear spline that uses a few unevenly spaced knots gives a fast and reasonably accurate method for computing $\Delta = \partial P/\partial S$.

The Black-Scholes price involves the parameters $r$ and $\sigma$, which need to be estimated. The yield of a short-maturity Treasury bill is usually used for $r$. Although in the GBM model for asset prices which are observed at fixed intervals of time (e.g. daily), one can estimate $\sigma$ by the standard deviation of historical (daily) asset returns, which are i.i.d. normal under the GBM model for asset prices, there are issues due to departures from this model (e.g., $\sigma$ can change over time and asset returns are markedly non-normal) and due to violations of the Black-Scholes assumptions in the financial market (e.g., there are actually transaction costs and limits on short selling). Section 13.4 and Chapter 16 of Hull [21] discuss how the parameter $\sigma$ in the Black-Scholes option price is treated in current practice. In the next section we describe an alternative approach that addresses the discrepancy between the Black-Scholes-Merton theory and time series data on American options and the underlying stock prices.

## 4. A combined substantive-empirical approach

In this section we describe an approach to time series modeling that contains both substantiative and empirical components. We first came up with this approach when we studied valuation of American options. Its basic idea is to use empirical modeling to address the gap between the actual prices in the American options market and the option prices given by the Black-Scholes-Merton theory in Example 2, as explained below.

**Example 3.** For European options, instead of using the basis function of Hutchinson, Lo and Poggio [22], an alternative approach is to express the option price as $c + Ke^{-rt^*}f^*(S/K, t^*)$, where $c$ is the Black-Scholes price (2.1) because the Black-Scholes formula has proved to be quite successful in explaining empirical data. This is tantamount to including $c(t, S)$ as one of the basis functions (with prescribed weight 1) to come up with a more parsimonious approximation to the actual option price.

The usefulness of this idea is even more apparent in the case of American options. Focusing on puts for definiteness, the decomposition formula (3.4) expresses an American put option price as the sum of a European put price $p$ and the early exercise premium which is typically small relative to $p$. This suggests that $p$ should be included as one of the basis functions (with prescribed weight 1). Lai and Wong [29] propose to use additive regression splines after the change of variables $u = -\sigma^2(T-t)$ and $z = \log(S/K)$. Specifically, for small $T-t$ (say within 5 trading days prior to expiration, i.e. $T-t \le 5/253$ under the assumption of 253 trading days per year), we approximate $P$ by $p$. For $T-t > 5/253$ (or equivalently, $u < -5\sigma^2/253$), we approximate $P$ by

$$\begin{aligned}
(4.1) \quad P = p &+ Ke^{\rho u}\{\alpha + \alpha_1 u + \sum_{j=1}^{J_u} \alpha_{1+j}(u - u^{(j)})_+ \\
&+ \beta_1 z + \beta_2 z^2 + \sum_{j=1}^{J_z} \beta_{2+j}(z - z^{(j)})^2_+ + \gamma_1 w + \gamma_2 w^2 \\
&+ \sum_{j=1}^{J_w} \gamma_{2+j}(w - w^{(j)})^2_+\},
\end{aligned}$$

where $\rho = r/\sigma^2$ as in (3.5), $\alpha$, $\alpha_j$, $\beta_j$ and $\gamma_j$ are regression parameters to be estimated by least squares from the training sample and

$$(4.2) \quad w = |u|^{-1/2}\{z - (\rho - \theta\rho - 1/2)u\} \quad (\theta = d/r)$$

is an "interaction" variable derived from $z$ and $u$. The motivation behind the centering term $(\rho - \theta\rho - 1/2)u$ comes from (3.5) that transforms GBM into Brownian motion, whereas that behind the normalization $|u|^{-1/2}$ comes from (3.4) and the closely related $d_1(x, y, v)$ in (2.2). The knots $u^{(j)}$ (respectively $z^{(j)}$ or $w^{(j)}$) of the linear (respectively quadratic) spline in (4.1) are the $100j/J_u$ (respectively $100j/J_z$ and $100j/J_w$)-th percentiles of $\{u_1, \ldots, u_n\}$ (respectively $\{z_1, \ldots, z_n\}$ or $\{w_1, \ldots, w_n\}$). The choice of $J_u$, $J_z$ and $J_w$ is over all possible integers between 1 and 10 to minimize the generalized cross validation (GCV) criterion, which can be expressed in the following form (cf. [19, 46]):

$$\text{GCV}(J_u, J_z, J_w) = \sum_{i=1}^n (P_i - \hat{P}_i)^2 \Big/ \left\{n\left(1 - \frac{J_u + J_z + J_w + 6}{n}\right)^2\right\},$$

where the $P_i$ are the observed American option prices in the past $n$ periods, and the $\hat{P}_i$ are the corresponding fitted values given by (4.1) in which the regression coefficients are estimated by least squares.

In the preceding we have assumed prescribed constants $\gamma$ and $\sigma$ as in the Black-Scholes model; these parameters appear in (4.1) via the change of variables (3.5). In practice $\sigma$ is unknown and may also vary with time. We can replace it in (4.1) by the standard deviation $\widehat{\sigma}_t$ of the most recent asset returns say, during the past 60 trading days prior to $t$ as in [22], p. 881. Moreover, the risk-free rate $r$ may also change with time, and can be replaced by the yield $\widehat{r}_t$ of a short-maturity Treasury bill on the close of the month before $t$. The same remark also applies to the dividend rate.

The simulation study in Lai and Wong [29] shows the advantages of this combined substantive-empirical approach. Not only is $P$ well approximated by $\widehat{P}$, especially over intervals of $S/K$ values that occur frequently in the sample, $\widehat{\Delta} - \Delta$ also reveals a pattern similar to that of $\widehat{P} - P$. Besides $\xi_{\widehat{P}} = E\{e^{-r\tau}|V_{\widehat{P}}(\tau)|\}$, where $\tau$ is the time of exercise and $V_{\widehat{P}}(t)$ is the value of the replicating portfolio at time $t$ that rebalances (according to the pricing formula $\hat{P}$) between the risky and riskless assets ([22], p. 868-869), Lai and Wong [29] also consider the measure

$$(4.3) \qquad \kappa_{\widehat{P}} = E\left\{\int_0^\tau (S_t/K)^2(\Delta(t) - \hat{\Delta}(t))^2 dt\right\},$$

where $\hat{\Delta} = \partial \hat{P}/\partial S$. In practice, continuous rebalancing is not possible. If rebalancing is done only daily, then $(S/K)^2(\Delta_A - \hat{\Delta})^2$ in (4.3) is replaced by a step function that stays constant on intervals of width $1/253$. Because of the adaptive nature of the methodology, the proposed approach of Lai and Wong [29] is much more robust to the misspecification error than the Black-Scholes formula in terms of both measures. Lai and Lim [26] carried out an empirical study of this approach and made use of its semiparametric pricing formula and (4.3) to come up with a modified Black-Scholes theory and optimal delta hedging in the presence of transaction costs.

## 5. Application to the 1821-1934 Canadian lynx data

The Canadian Lynx data set consists of the annual record of the numbers of the Canadian lynx trapped in the Mackenzie River district of the North-west Canada for the period 1821-1934 inclusively. Let $X_t$ be $\log_{10}$(number recorded as trapped in year $1820+t$) ($t = 1, \ldots, 114$). Figure 1 shows the time series plot of $X_t$. According to Tong [44], Moran [39] performed the first time series analysis on these data by fitting an AR(2) model to $X_t$; moreover, the log transformation is used because it (i) makes the marginal distribution of $X_t$ more symmetric about its mean and (ii) reduces the approximation error in assuming the number of lynx to be proportional to the population. In view of the substantial non-linearity of $E[X_t|X_{t-3}]$ found in the scatterplot of $X_t$ versus $X_{t-3}$, Tong([44], p.361) critiques Moran's analysis and its enhancements by Campbell and Walker [6], who added a harmonic component to the AR(2) model, and by Tong [43], who used the AIC to select the order $p = 11$ for AR(p) models, as "uncritical acceptance of linearity" in $X_t$. He uses a self-excited threshold autoregressive model (SETAR) of the form

$$(5.1) \qquad X_t - X_{t-1} = \begin{cases} 0.62 + 0.25 X_{t-1} - 0.43 X_{t-2} + \varepsilon_t & \text{if } X_{t-2} \leq 3.25 \\ -(1.24 X_{t-2} - 2.25) + 0.52 X_{t-1} + \varepsilon_t & \text{if } X_{t-2} > 3.25 \end{cases}$$

to fit these data, similar to Tong and Lim ([45], Section 9). The growth rate $X_t - X_{t-1}$ in the first regime (i.e., $X_{t-2} \leq 3.25$) tends to be positive but small, which corresponds to a slow population growth. In the second regime (i.e., $X_{t-2} > 3.25$), $X_t - X_{t-1}$ tends to be negative, corresponding to a decrease in population size.

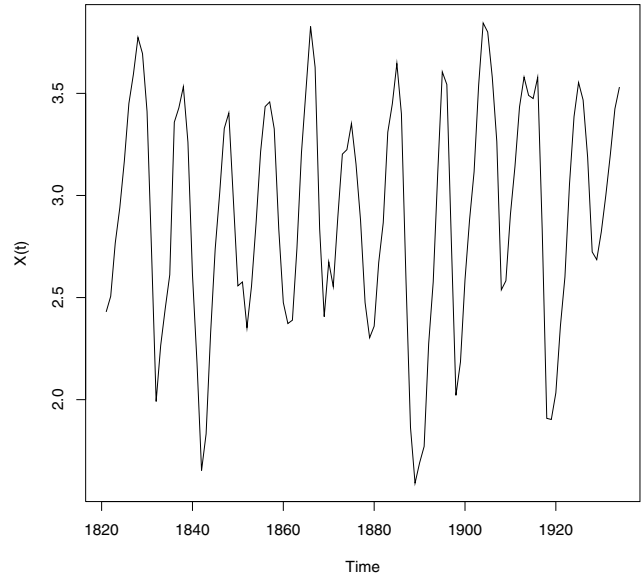

FIG 1. *Time series plot of $\log_{10}$ of the Canadian lynx series.*

Tong ([44], p. 377) interprets the fitted model as an "energy balance" between the population expansion and the population contraction, yielding a stable limit cycle with a 9-year period which is in good agreement with the observed asymmetric cycles. Motivated by Van der Pol's equation, Haggan and Ozaki [18] proposed to fit another nonlinear time series model, namely, the exponential autoregressive model

$$(5.2) \qquad X_t - \mu = \sum_{j=1}^{11}(\phi_j + \pi_j e^{-\gamma(X_{t-j}-\mu)^2})(X_{t-j}-\mu) + \varepsilon_t,$$

which gives a limit cycle of period 9.45 years. Lim [35] compares the prediction performance of these and other parametric models and concludes that Tong's SETAR model ranks the best among them.

Taking a more nonparametric approach, Fan and Yao [14] use a functional − coefficient autoregressive model to fit the observed $X_t$ series and compare its prediction with that of threshold autoregression. Specifically, they fit the FAR(2,2) model

$$(5.3) \qquad X_t = a_1(X_{t-2})X_{t-1} + a_2(X_{t-2})X_{t-2} + \sigma\varepsilon_t$$

to the first 102 observations, reserving the last 12 observations to evaluate the prediction. The $a_1(\cdot)$ and $a_2(\cdot)$ in (5.3) are unknown functions which are estimated by using locally linear smoothers. Fan and Yao ([14], p. 327) plot the estimates $\hat{a}_1(\cdot)$ and $\hat{a}_2(\cdot)$, which are approximately constant for $X_{t-2} < 2.7$ with $\hat{a}_1(X_{t-2}) \approx 1.3$ and $\hat{a}_2(X_{t-2}) \approx -0.2$, and which are approximately linear for $X_{t-2} \geq 2.7$. For comparison, Fan and Yao [14] also fit the following SETAR(2) model to the same set of data:

$$(5.4) \qquad \widehat{X}_t = \begin{cases} 0.424 + 1.255 X_{t-1} - 0.348 X_{t-2}, & X_{t-2} \leq 2.981, \\ 1.882 + 1.516 X_{t-1} - 1.126 X_{t-2}, & X_{t-2} > 2.981. \end{cases}$$

Because of the close resemblance of the fitted SETAR(2) and FAR(2,2), they share certain ecological interpretations. In particular, the difference of the fitted coefficients in each regime can be explained by using the phase dependence and the density dependence in the predator-prey structure. The phase dependence refers to the difference in the behavior of preys (snowshoe hare) and predators (lynx) in hunting and escaping at the decreasing and increasing phases of population dynamics, while the density dependence is the relationship between the reproduction rates of the animals and their abundance. More discussion on these ecological interpretations can be found in [42].

To evaluate the predictions of FAR (2,2), Fan and Yao ([14], p. 324) use the one-step ahead forecast (denoted by $\widehat{X}_t$) and the iterative two-step-ahead forecast (denoted by $\tilde{X}_t$), which are defined by

$$\widehat{X}_t := \hat{a}_1(X_{t-2})X_{t-1} + \hat{a}_2(X_{t-2})X_{t-2}, \quad \tilde{X}_t := \hat{a}_1(X_{t-2})\widehat{X}_{t-1} + \hat{a}_2(X_{t-2})X_{t-2}.$$

The predictions of SETAR(2) are similarly defined. The out-sample prediction absolute errors ($|\hat{X}_t - X_t|$ and $|\tilde{X}_t - X_t|$) of the last 12 observations are reported in Table 1. Based on the average of these 12 absolute prediction errors (AAPE), FAR(2,2) performs slightly better than SETAR(2). Other nonparametric time series models for the Canadian lynx data include the projection pursuit regression (PPR) model fitted by Lin and Pourahmadi [36] who found that SETAR outperforms PPR in terms of one-step-ahead forecasts, and neural network models which Kajitani, Hipel and McLeod [25] found to be "just as good or better than SETAR models for one-step out-of-sample forecasting of the lynx data."

A substantive approach is adopted by Royama ([41], Chapter 5). Instead of building the statistical model first and using ecology to interpret the fitted model later, Royama starts with ecological mechanisms and population dynamics. Letting $R_t = X_{t+1} - X_t$ denote the log reproductive rate from year $t$ to $t+1$, he considers nonlinear dynamics of the form $R_t = f(X_t, \ldots, X_{t-h+1}) + u_t$, where $u_t$ is a zero-mean random disturbance, and emphasizes that "our ultimate goal is to determine the reproduction surface $f$ and to find an appropriate model which reasonably approximates to it," with $f$ satisfying the following two conditions in view of ecological considerations: There exists $X^*$ such that $f(X^*, \ldots, X^*) = 0$, and $R_t$ has to be bounded above because "no animal can produce infinite number of offspring" (see [41], p. 50, 154, 178). In Chapter 4 of [42], Royama introduces the (first-order) logistic model of $f(X_t) = r_m - \exp\{-a_0 - a_1 X_{t-1}\}$ to incorporate competition over an available resource. Here $r_m$ is the maximum biologically realizable reproduction rate, i.e. $R_t \leq r_m$ for all $t$; see [42], Section 4.2.5. An implicit assumption of the model is that the resource being depleted during a time step will be recovered to the same level by the onset of the next time step. This assumption can be relaxed if a linear combination of $X_{t-j}(j = 1, \ldots, h)$ with $h > 1$ is used in the exponential term of $f$, yielding a higher-order logistic model; see [41], p. 153.

TABLE 1
*Absolute prediction errors of one-step-ahead (1 yr) and iterative two-step-ahead (2 yr) forecasts and their 12-year average (AAPE).*

|      |       | Model (5.3) FAR(2,2) | | Model (5.4) SETAR(2) | | Model (5.6) Logistic | | Model (5.8a) Logistic-MARS | |
|------|-------|------|------|------|------|------|------|------|------|
| Year | $X_t$ | 1 yr | 2 yr | 1 yr | 2 yr | 1 yr | 2 yr | 1 yr | 2 yr |
| 1923 | 3.054 | 0.157 | 0.156 | 0.187 | 0.090 | 0.178 | 0.075 | 0.188 | 0.082 |
| 1924 | 3.386 | 0.012 | 0.227 | 0.035 | 0.269 | 0.077 | 0.281 | 0.057 | 0.286 |
| 1925 | 3.553 | 0.021 | 0.035 | 0.014 | 0.038 | 0.057 | 0.153 | 0.073 | 0.120 |
| 1926 | 3.468 | 0.008 | 0.037 | 0.022 | 0.000 | 0.012 | 0.077 | 0.023 | 0.140 |
| 1927 | 3.187 | 0.085 | 0.101 | 0.059 | 0.092 | 0.020 | 0.018 | 0.122 | 0.168 |
| 1928 | 2.723 | 0.055 | 0.086 | 0.075 | 0.015 | 0.128 | 0.098 | 0.002 | 0.159 |
| 1929 | 2.686 | 0.135 | 0.061 | 0.273 | 0.160 | 0.179 | 0.004 | 0.009 | 0.012 |
| 1930 | 2.821 | 0.016 | 0.150 | 0.026 | 0.316 | 0.004 | 0.216 | 0.010 | 0.001 |
| 1931 | 3.000 | 0.017 | 0.037 | 0.030 | 0.062 | 0.005 | 0.010 | 0.013 | 0.025 |
| 1932 | 3.201 | 0.007 | 0.014 | 0.060 | 0.043 | 0.048 | 0.042 | 0.021 | 0.005 |
| 1933 | 3.424 | 0.089 | 0.098 | 0.076 | 0.067 | 0.124 | 0.184 | 0.066 | 0.091 |
| 1934 | 3.531 | 0.053 | 0.175 | 0.072 | 0.187 | 0.083 | 0.245 | 0.011 | 0.087 |
| AAPE |       | 0.055 | 0.095 | 0.073 | 0.112 | 0.075 | 0.117 | 0.050 | 0.098 |

Chapter 5 of Royama [41] examines the autocorrelation function and the partial autocorrelation function of the Canadian lynx series and concludes that $h$ should be set to 2, which corresponds to the model

$$(5.5) \qquad X_t - X_{t-1} = r_m - \exp\{-a_0 - a_1 X_{t-1} - a_2 X_{t-2}\} + u_{t-1},$$

where $r_m, a_0, a_1$ and $a_2$ are unknown parameters that need to be estimated; see [41], p. 190-191. From the scatterplot of $R_{t-1} = X_t - X_{t-1}$ versus $X_{t-2}$, Royama guesses $r_m \approx 0.6$ and $X^* \approx 3$. He uses this together with trial and error to obtain the estimate $(\widehat{r}_m, \widehat{a}_0, \widehat{a}_1, \widehat{a}_2) = (0.597, 2.526, 0.838, -1.508)$, but finds that the asymmetric cycle of the fitted model does not match the observed cycle from the data well. Moreover, the fitted autocorrelation function decays too fast in comparison with the sample autocorrelation function.

Instead of his *ad hoc* estimates, we can use nonlinear least squares, initialized at his estimates, to estimate the parameters of (5.5), yielding

$$(5.6) \qquad X_t - X_{t-1} = 0.460 - \exp\{-3.887 - 0.662 X_{t-1} + 1.663 X_{t-2}\} + u_{t-1},$$

which implies that the maximum logarithmic reproduction rate is 0.460, i.e., the population can grow at most $10^{0.46} = 2.884$ times per year. Figure 2, top left corner, shows a negative contour of the response surface of the fitted model (5.6). This implies that the population size can drop sharply in the region $X_{t-2} > 3.5$ and $X_{t-1} < 2.5$, leading to extinction in the upper left part of this region. Whereas (5.6) does not rule out the possibility of $X_t$ diverging to $-\infty$, extinction occurs as soon as $X_t$ falls below 0 (or equivalently, the population size $10^{X_t}$ falls below 1).

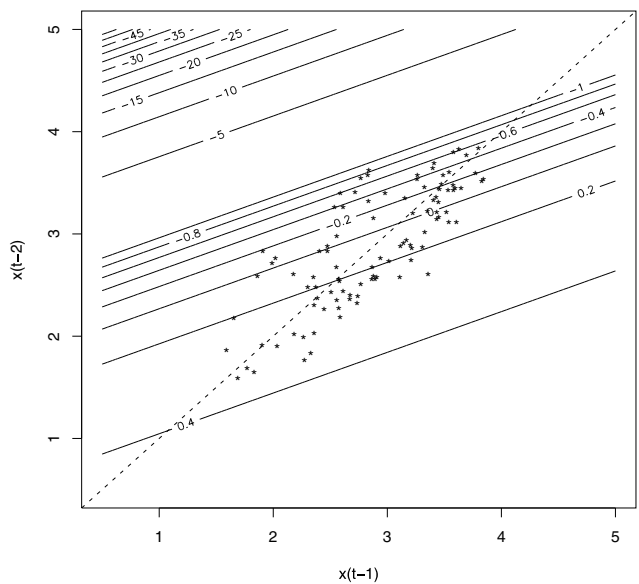

FIG 2. *Contour plot of $\hat{R}_{t-1} = \widehat{X_t - X_{t-1}}$ of the logistic model (5.6). The observations are marked by $*$. The dotted line is $X_{t-2} = X_{t-1}$. The intersection of this line and the contour numbered 0 gives the equilibrium $X^*$.*

Note that one can also derive bounds on the logarithmic reproduction rates from the empirical approach. Figure 3 is the plot of the limit cycle generated by the skeleton of the fitted model (5.4). The limit cycle is of period 8 years. The maximum and the minimum logarithmic reproduction rates, attained at years 1 and 5 in Figure 3, are 0.212 and -0.269, respectively. That is, the population grows at most $10^{0.212} = 1.629$ times per year and diminishes by at most a factor of $10^{-0.269}=0.538$ per year. Moreover, the limit cycle of (5.4) implies an infinite loop of expansion and contraction and rules out the possibility of extinction. These are consequences of adopting an empirical approach because the data are distributed along the main diagonal of Figure 2, but not its top left corner nor its lower right corner. In order to deduce the behavior of the reproduction rates in these regions, mechanistic modeling is essential. On the other hand, the empirical approach uses the observed data better and gives more accurate forecasts. Table 1 compares the prediction performance of FAR(2,2) and SETAR(2) with that of the logistic model (5.5). The fitted logistic model provides the worst AAPE of one-step-ahead and iterative two-step-ahead forecasts. Moreover, instead of characterizing the equilibrium with limit cycles, the logistic model only gives two equilibrium points, with one corresponding to extinction and the other equal to $X^* = \{a_0 + \log(r_m)\}/(a_1 + a_2) = 3.107$ (the intersection of the line $X_{t-1} = X_{t-2}$ and the contour of $f = 0$ in Figure 2.)

We next apply the combined substantive-empirical approach of Section 4 to these data, using the substantive model (5.5) to provide one of the basis functions in the

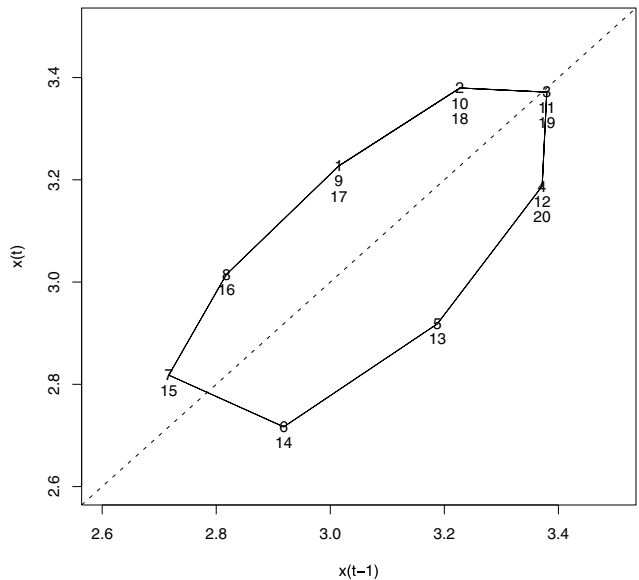

FIG 3. *Limit cycle of the skeleton of the SETAR(2) model (5.4). The dotted line is $X_t = X_{t-1}$.*

semiparametric model

$$\text{(5.7)} \quad \begin{aligned} X_t - X_{t-1} &= r_m - \exp\{-a_0 - a_1 X_{t-1} - a_2 X_{t-2}\} \\ &\quad + g(X_{t-1}, X_{t-2}) I\{(X_{t-1}, X_{t-2}) \in S\} + u_{t-1}, \end{aligned}$$

where $g$ is an unknown function and $S$ is a region containing the observations that will be specified later. Moreover, the difference equation (5.7) has the boundary constraint $X_{t-1} + r_m - \exp\{-a_0 - a_1 X_{t-1} - a_2 X_{t-2}\} + g(X_{t-1}, X_{t-2}) I\{(X_{t-1}, X_{t-2}) \in S\} \geq 0$. The lynx population becomes extinct as soon as this boundary condition is violated. Model (5.7) can be fitted by using the backfitting algorithm. Specifically, model (5.5) is estimated first and then the residuals are used as the response variable in nonparametric regression on the predictor variable $(X_{t-1}, X_{t-2})$. The difference between the observed $X_t - X_{t-1}$ and the fitted $g$ is then used as the response variable in (5.5), whose parameters can be estimated by nonlinear least squares. The algorithm of multivariate adaptive regression splines (MARS) developed by Friedman (1991) is used for estimating $g$ for the first step in each iteration of the above backfitting procedure (the function "mars" in the package of "mda" in R can be used). This kind of iteration sheme has been used in fitting *partly linear* models, where the parametric component is a linear regression model and the nonparametric component is often fitted by using kernel regression; see [8, 13, 20]. The fitted response surface is

$$\text{(5.8a)} \quad \begin{aligned} X_t - X_{t-1} &= 1.319 - \exp\{-0.224 - 0.205 X_{t-1} + 0.343 X_{t-2}\} \\ &\quad + \hat{g}(X_{t-1}, X_{t-2}) I\{(X_{t-1}, X_{t-2}) \in S\} + u_{t-1}, \end{aligned}$$

$$\text{(5.8b)} \quad \begin{aligned} \hat{g}(X_{t-1}, X_{t-2}) &= 2.294 (X_{t-1} - 3.224)_+ (X_{t-2} - 2.864)_+ \\ &\quad - 1.572 (X_{t-1} - 3.202)_+ - 0.851 (X_{t-2} - 3.202)_+. \end{aligned}$$

We evaluate this fitted model by using the out-sample prediction criterion. Table 1 shows that (5.8a) gives the smallest AAPE for one-step-ahead forecasts among all models considered, and that the AAPE for iterative two-step-ahead forecasts of (5.8a) is comparable to the smallest one provided by FAR(2,2). The region $S$ in (5.8a) is chosen to be the oblique rectangle whose edges are defined by the sample means $\pm 3$ standard deviations of the principal components of the bivariate sample of $(X_{t-1}, X_{t-2})$; see Figure 4 which shows that this region contains not only the in-sample data but also the out-sample data. Figure 5 gives the contour plot of the fitted model (5.8a). The logarithmic growth rate at its top left corner is about $-2$, which shows a strong possibility of extinction even though the magnitude is less drastic than that in Figure 2 for (5.6). The inclusion of tensor products of univariate splines in (5.8a) would have produced positive probability limits of $X_t$ diverging to $\infty$ or to $-\infty$ if $(X_{t-1}, X_{t-2})$ had not been confined to a compact region. On the other hand, with an absorbing barrier at 0 and with (5.8b) only applicable inside the compact set $S$, Markov chains of the type (5.8a) not only have stationary distributions but are also geometrically ergodic under mild assumptions on the random disturbances $u_t$ (e.g., to ensure irreducibity); see [39].

## 6. Conclusion

In his concluding remarks, Cox [11] noted that for successful use of statistical models in particular applications, "large elements of subject-matter judgment and technical statistical expertise are usually essential. Indeed, it is precisely the need for this combination that makes our subject such an interesting and demanding one." We have followed up on his remarks here with a combined subject-matter and statistical modeling approach to time series analysis, which we illustrate for the "particular applications" of option pricing and population dynamics of the Canadian lynx. In particular, for the Canadian lynx data, we have shown how statistical modeling for data-rich regions of $(X_{t-1}, X_{t-2})$ can be combined effectively with "subject-matter judgment" which is the only reliable guide for sparse-data regions.

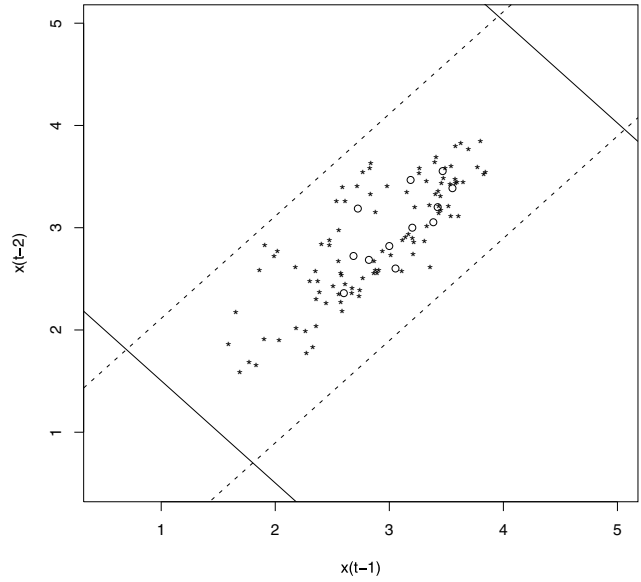

FIG 4. *The oblique rectangle $S$ formed by $\pm 3$ standard deviations away from the sample means of the principal components of $(X_{t-1}, X_{t-2})$. The in-sample and out-sample observations are marked by $*$ and o, respectively.*

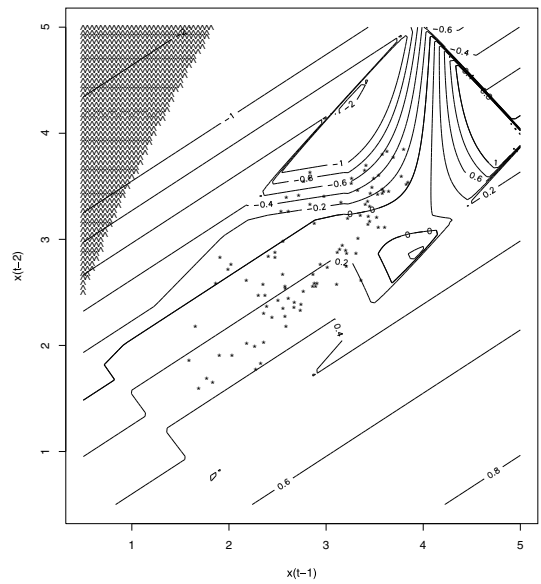

FIG 5. *Contour plot of $\hat{R}_{t-1} = X_t - \widehat{X_{t-1}}$ of the logistic-MARS model (5.7). The observations are marked by $*$. The shaded region corresponds to extinction.*

**Acknowledgments**

Lai's research was supported by the National Science Foundation grant DMS-0305749. Wong's research was supported by the Research Grants Council of Hong Kong under grant CUHK6158/02E.

**References**

[1] AITSAHLIA, F. AND LAI, T. L. (2001). Exercise boundaries and efficient approximations to American option prices and hedge parameters. *J. Comput. Finance* **4** 85–103.
[2] BARRON, A. R. (1993). Universal approximation bounds for superpositions of a sigmoid function. *IEEE Trans. Information Theory* **39** 930–945. MR1237720
[3] BLACK, F. AND SCHOLES, M. (1973). The pricing of options and corporate liabilities. *J. Political Economy* **81** 637–659.
[4] BROADIE, M. AND DETEMPLE, J. (1996). American option valuation: New bounds, approximations, and a comparison of existing methods. *Rev. Financial Studies* **9** 1121–1250.
[5] BROADIE, M., DETEMPLE, J., GHYSELS, E. AND TORRES, O. (2000). Nonparametric estimation of American options' exercise boundaries and call prices. *J. Econ. Dynamics & Control* **24** 1829–1857. MR1784575
[6] CAMPBELL, M. J. AND WALKER, A.M. (1977). A survey of statistical work on the McKenzie River series of annual Canadian lynx trappings for the years 1821-1934, and a new analysis. *J. Roy. Statist. Soc. Ser. A* **140** 411–431.
[7] CARR, P., JARROW, R. AND MYNENI, R. (1992). Alternative characterizations of American put options. *Math. Finance* **2** 87–106. MR1143390
[8] CHEN, H. (1988). Convergence rates for parametric components in a partly linear model. *Ann. Statist.* **16** 136–146. MR0924861
[9] CHEN, R. AND TSAY, R. S. (1993). Functional-coefficient autoregressive models. *J. Amer. Statist. Assoc.* **88** 298–308. MR1212492
[10] CHEN, R. AND TSAY, R. S. (1993). Nonlinear additive ARX models. *J. Amer. Statist. Assoc.* **88** 955–967.
[11] COX, D. R. (1990). Role of models in statistical analysis. *Statist. Sci.* **5** 169–174. MR1062575
[12] COX, J., ROSS, S. AND RUBINSTEIN, M. (1979). Option pricing: A simplified approach. *J. Financial Econ.* **7** 229–263.
[13] ENGLE, R. F., GRANGER, C. W. J., RICE, J. AND WEISS, A. (1986). Semiparametric estimates of the relation between weather and electricity sales. *J. Amer. Statist. Assoc.* **81** 310–320.
[14] FAN, J. AND YAO, Q. (2003). *Nonlinear Time Series.* Springer-Verlag, New York. MR1964455
[15] FRIEDMAN, J. H. (1991). Multivariate adaptive regression splines. *Ann. Statist.* **19** 1–142. MR1091842

[16] GOODWIN, G. C., RAMADGE, P. J. AND CAINES, P. E. (1981). Discrete time stochastic adaptive control. *SIAM J. Control Optim.* **19** 829–853. MR0634955

[17] GUO, L. AND CHEN, H. F. (1991). The Åström-Wittenmark self-tuning regulator revisited and ELS-based adaptive trackers. *IEEE Trans. Automat. Contr.* **36** 802–812. MR1109818

[18] HAGGAN, V. AND OZAKI, T. (1981). Modelling non-linear random vibrations using an amplitude-dependent autoregressive time series model. *Biometrika* **68** 189–196. MR0614955

[19] HASTIE, T. J. AND TIBSHIRANI, R. J. (1990). *Generalized Additive Models.* Chapman & Hall, London. MR1082147

[20] HECKMAN, N. E. (1988). Minimax estimates in a semiparametric model. *J. Amer. Statist. Assoc.* **83** 1090–1096. MR0997587

[21] HULL, J. C. (2006). *Options, Futures and Other Derivatives*, 6th edn. Pearson Prentice Hall, Upper Saddle River, NJ.

[22] HUTCHINSON, J. M., LO, A. W. AND POGGIO, T. (1994). A nonparameric approach to pricing and hedging derivative securities via learning networks. *J. Finance* **49** 851–889.

[23] JACKA, S. D. (1991). Optimal stopping and the American put. *Math. Finance* **1** 1–14.

[24] JU, N. (1998). Pricing an American option by approximating its early exercise boundary as a multipiece exponential function. *Rev. Financial Studies* **11** 627–646.

[25] KAJITANI, Y., HIPELM, K. W. AND MCLEOD, A. I. (2005). Forecasting nonlinear time series with feed-forward neural networks: A case study of Canadian Lynx data. *J. Forecasting* **24** 105–117. MR2148983

[26] LAI, T. L. AND LIM, T. W. (2006). A new approach to pricing and hedging options with transaction costs. Tech. Report, Dept. Statistics, Stanford Univ.

[27] LAI, T. L. AND WEI, C. Z. (1987). Asymptotically efficient self-tuning regulators. *SIAM J. Control Optim.* **25** 466–481. MR0877072

[28] LAI, T. L. AND WONG, S. P. (2001). Stochastic neural networks with applications to nonlinear time series. *J. Amer. Statist. Assoc.* **96** 968–981. MR1946365

[29] LAI, T. L. AND WONG, S. P. (2004). Valuation of American options via basis functions. *IEEE Trans. Automat. Contr.* **49** 374–385. MR2062250

[30] LAI, T. L. AND YING, Z. (1991). Parallel recursive algorithms in asymptotically efficient adaptive control of linear stochastic systems. *SIAM J. Control Optim.* **29** 1091–1127. MR1110088

[31] LEHMANN, E. L. (1990). Model specification: The views of Fisher and Neyman, and later developments. *Statist. Sci.* **5** 160–168. MR1062574

[32] LEWIS, P. A. W. AND RAY, B. K. (1993). Nonlinear modeling of multivariate and categorical time series using multivariate adaptive regression splines. In *Dimension Estimation and Models* (H. Tong, ed). World Sci. Publishing, River Edge, NJ, pp. 136–169. MR1307658

[33] LEWIS, P. A. W. AND RAY, B. K. (2002). Nonlinear modelling of periodic threshold autoregressions using TSMARS. *J. Time Ser. Anal.* **23** 459–471. MR1910892

[34] LEWIS, P. A. W. AND STEVENS, J. G. (1991). Nonlinear modeling of time series using multivariate adaptive regression splines (MARS). *J. Amer. Statist. Assoc.* **86** 864–877.

[35] LIM, K. S. (1987). A comparative study of various univariate time series models for Canadian lynx data. *J. Time Ser. Anal.* **8** 161–176.

[36] LIN, T. C. AND POURAHMADI, M. (1998). Nonparametric and nonlinear mod-

els and data mining in time series: a case-study on the Canadian lynx data. *Appl. Statist.* **47** 187–201.
[37] MERTON, R. C. (1973). Theory of rational option pricing. *Bell J. Econ. & Management Sci.* **4** 141–181. MR0496534
[38] MEYN, S. P. AND TWEEDIE, R. L. (1993). *Markov Chains and Stochastic Stability.* Springer-Verlag, New York. MR1287609
[39] MORAN, P. A. P. (1953). The statistical analysis of the Canadian lynx cycle, I: Structure and prediction. *Austral. J. Zoology* **1** 163–173.
[40] ROSS, S. A. (1987). Finance. In *The New Palgrave: A Dictionary of Economics* (J. Eatwell, M. Milgate and P. Newman, eds.), Vol. 2. Stockton Press, New York, pp. 322–336.
[41] ROYAMA, T. (1992). *Analytical Population Dynamics.* Chapman & Hall, London.
[42] STENSETH, N. C., CHAN, K. S., TONG, H., BOONSTRA, R., BOUTIN, S., KREBS, C. J., POST, E., O'DONOGHUE, M., YOCCOZ, N. G., FORCHHAMMER, M. C. AND HURRELL, J. W. (1998). From patterns to processes: Phase and density dependencies in the Canadian lynx cycle. *Proc. Natl. Acad. Sci. USA* **95** 15430–15435.
[43] TONG, H. (1977). Some comments on the Canadian lynx data. *J. Roy. Statist. Soc. Ser. A* **140** 432–435.
[44] TONG, H. (1990). *Nonlinear Time Series.* Oxford University Press, Oxford. MR1079320
[45] TONG H. AND LIM, K. S. (1980). Threshold autoregression, limit cycles and cyclical data (with Discussion). *J. Roy. Statist. Soc. Ser. B* **42** 245–292.
[46] WAHBA, G. (1990). *Spline Models for Observational Data.* SIAM Press, Philadelphia. MR1045442
[47] WEIGEND, A. AND GERSHENFELD, N. (1993). *Time Series Prediction: Forecasting the Future and Understanding the Past.* Addison-Wesley, Reading, MA.
[48] WEIGEND, A., RUMELHART, D. AND HUBERMAN, B. (1991). Predicting Sunspots and Exchange Rates with Connectionist Networks. In *Nonlinear Modeling and Forecasting* (Casdagli, M. and Eubank, S., eds.). Addison Wesley, Redwood City, CA, 395–432.